\documentclass[11pt]{amsart}

\newtheorem{thm}{Theorem}[section]
\newtheorem{cor}[thm]{Corollary}
\newtheorem{lem}[thm]{Lemma}
\newtheorem{proposition}[thm]{Proposition}

\theoremstyle{definition}

\newtheorem{rem}[thm]{Remark}

\numberwithin{equation}{section}

\usepackage[active]{srcltx} 

\def\J#1#2#3{ \left\{ #1,#2,#3 \right\} }

\def\NN{{\mathbb{N}}}
\def\11{{1}}
\def\CC{{\mathbb{C}}}

\usepackage{enumerate}
\usepackage{amsmath}
\usepackage{amssymb}
\usepackage{hyperref}

\begin{document}


\baselineskip=17pt



\title[Automatic continuity of biorthogonality preservers]{Automatic continuity of biorthogonality
preservers between weakly  compact JB$^*$-triples and atomic JBW$^*$-triples}

\author[Burgos]{Mar{\' i}a Burgos}
\email{mariaburgos@ugr.es}
\address{Departamento de An{\'a}lisis Matem{\'a}tico, Facultad de
Ciencias, Universidad de Granada, 18071 Granada, Spain.}

\author[Garc{\' e}s]{Jorge J. Garc{\' e}s}
\email{jgarces@ugr.es}
\address{Departamento de An{\'a}lisis Matem{\'a}tico, Facultad de
Ciencias, Universidad de Granada, 18071 Granada, Spain.}

\author[Peralta]{Antonio M. Peralta}
\email{aperalta@ugr.es}
\address{Departamento de An{\'a}lisis Matem{\'a}tico, Facultad de
Ciencias, Universidad de Granada, 18071 Granada, Spain.}

\thanks{Published at Studia Mathematica \url{https://doi.org/https://doi.org/10.4064/sm204-2-1}. This manuscript version is made available under the CC-BY-NC-ND 4.0 license \url{https://creativecommons.org/licenses/by-nc-nd/4.0/}}

\date{}

\begin{abstract}
We prove that every biorthogonality preserving
linear \hyphenation{sur-jection} surjection from a weakly compact JB$^*$triple containing no
infinite dimensional rank-one summands onto another JB$^*$-triple
is automatically continuous. We also show that every
biorthogonality preserving linear surjection between atomic
JBW$^*$triples containing no infinite dimensional rank-one
summands is automatically continuous. Consequently, two \hyphenation{ato-mic}atomic 
JBW$^*$-triples containing no rank-one summands are isomorphic if,
and only if, there exists a (non necessarily continuous)
biorthogonality preserving linear surjection between them.
\end{abstract}



\maketitle

\section{Introduction and preliminaries}

Studies on the automatic continuity of those linear surjections
between C$^*$- and von Neumann algebras preserving orthogonality
relations in both directions constitute the latest variant of a
problem initiated by W. Arendt in early eighties.\smallskip

We recall that two complex valued continuous functions $f$ an $g$
are said to be \emph{orthogonal} whenever they have disjoint
supports. A mapping $T$ between $C(K)$ spaces is called
orthogonality preserving if it maps orthogonal functions to
orthogonal functions. The main result established by Arendt proves
that every orthogonality preserving bounded linear mapping $T :
C(K) \to C(K)$ is of the form $$T(f) (t) = h(t) f(\varphi (t)) \
(f\in C(K), \ t\in K),$$ where $h\in C(K)$ and $\varphi: K \to K$
is a mapping which is continuous on $\{ t\in K : h(t) \neq
0\}$.\smallskip

The hypothesis of $T$ being continuous was relaxed by K. Jarosz in
\cite{Jar}. In fact, Jarosz obtained a complete description of
every orthogonality preserving (non-necessarily continuous) linear
mapping between $C(K)$-spaces. Among the consequences of his
description, it follows that an orthogonality preserving linear
surjection between $C(K)$-spaces is automatically
continuous.\smallskip

Two elements $a,b$ in a general C$^*$-algebra $A$ are said to be
\emph{orthogonal} (denoted by $a \perp b$) if $a b^* = b^* a=0$.
When $a= a^*$ and $b=b^*$, we have $a\perp b$ if and only if $a
b=0$. A mapping $T$ between two C$^*$-algebras $A$, $B$ is called
\emph{orthogonality preserving} if $T(a) \perp T(b)$ for every
$a\perp b$ in $A$. When $T(a) \perp T(b)$ in $B$ if and only if
$a\perp b$ in $A$, we say that $T$ is \emph{biorthogonality
preserving}. Under continuity assumptions, orthogonality
preserving bounded linear operators between C$^*$-algebras are
completely described in \cite[\S 4]{BurFerGarMarPe}. The just
quoted paper culminates the studies developed by W. Arendt
\cite{Arendt}, K. Jarosz \cite{Jar}, M. Wolff \cite{Wol}, and N.-Ch.
Wong \cite{Won}, among others, on bounded orthogonality
preserving linear maps between C$^*$-algebras.\smallskip

C$^*$-algebras belongs to a wider class of complex Banach spaces
in which orthogonality also makes sense. We refer to the class of
(complex) JB$^*$-triples (see \S 2 for definitions). Two elements
$a,b$ in a JB$^*$-triple $E$ are said to be \emph{orthogonal}
(denoted by $a \perp b$) if $L(a,b)=0,$ where $L(a,b) $ is the
linear operator in $E$ given by $L(a,b)x=\{a,b,x\}$. A linear
mapping $T: E \rightarrow F$ between two JB$^*$-triples is called
\emph{orthogonality preserving} if $T(x) \perp T(y)$ whenever $x
\perp y$. The mapping $T$ is \emph{biorthogonality preserving}
whenever the equivalence $ x \perp y \Leftrightarrow T(x) \perp
T(y)$ holds for all $x,y$ in $E$.\smallskip

Most of the novelties introduced in \cite{BurFerGarMarPe} consist
in studying orthogonality preserving bounded linear operators from
a C$^*$-algebra or a JB$^*$-algebra to a JB$^*$-triple to take
advantage of the techniques developed in JB$^*$-triple theory.
These techniques were successfully applied in the subsequent paper
\cite{BurFerGarPe} to obtain the description of such operators
(see \S 2 for a detailed explanation).\smallskip

Despite the vast literature which has appeared on orthogonality
preserving bounded\hyphenation{boun-ded} linear operators between
C$^*$-algebras and JB$^*$-triples, just a reduced number of papers
have considered the problem of automatic continuity of
biorthogonality preserving linear surjections between
C$^*$-algebras. Besides Jarosz \cite{Jar}, mentioned above,
M. A. Chebotar, W.-F. Ke, P.-H. Lee, and N.-C.
Wong proved, in \cite[Theorem 4.2]{ChebKeLeeWo}, that every
zero-products preserving linear bijection from a properly infinite
von Neumann algebra into a unital ring is a ring homomorphism followed
by left multiplication by the image of the identity.
J. Araujo and K. Jarosz showed that every linear
bijection between algebras $L(X),$ of continuous linear maps on a
Banach space $X$, which preserves zero-products in both directions
is automatically continuous and a multiple of an algebra
isomorphism \cite{ArauJar}. These authors also conjectured that
every linear bijection between two C$^*$-algebras preserving
zero-products in both directions is automatically continuous (see
\cite[Conjecture 1]{ArauJar}).\smallskip

The authors of this note proved in \cite{BurGarPe} that every
biorthogonality preserving linear surjection between two compact
C$^*$-algebras or between two von Neumann algebras is
automatically continuous. One of the consequences of this result
is a partial answer to \cite[Conjecture 1]{ArauJar}. Concretely,
every surjective and symmetric linear mapping between von Neumann
algebras (or compact C$^*$-algebras) which preserves zero-products
in both directions is continuous.\smallskip

In this paper we study the problem of automatic continuity of
biorthogonality preserving linear surjections between
JB$^*$-triples, extending some of the results obtained in
\cite{BurGarPe}. Section 2 contains the basic definitions and
results used in the paper. Section 3 is devoted to survey the
structure and properties of the (orthogonal) annihilator of a
subset $M$ in a JB$^*$-triple paying attention to the annihilator
of a single element. In section 4 we prove that every
biorthogonality preserving linear surjection from a weakly compact
JB$^*$-triple containing no infinite dimensional rank-one summands
to a JB$^*$-triple is automatically continuous. In Section 5 we
show that two atomic JB$^*$-triples containing no rank-one
summands are isomorphic if, and only if, there exists a
biorthogonality preserving linear surjection between them, a
result which follows from the automatic continuity of every
biorthogonality preserving linear surjection between atomic
JB$^*$-triples containing no infinite dimensional rank-one
summands.

\section{Notation and preliminaries}

Given Banach spaces $X$ and $Y$, $L(X,Y)$ will denote the space of
all bounded linear mappings from $X$ to $Y$. The symbol $L(X)$
will stand for the space $L(X,X)$. Throughout the paper the word
``operator'' will always mean bounded linear mapping. The dual
space of a Banach space $X$ is denoted by $X^*$.\smallskip

 JB$^*$-triples were introduced
by W. Kaup in \cite{Ka}. A JB$^*$-triple is a complex Banach space
$E$ together with a continuous triple product $\J ... : E\times
E\times E \to E,$ which is conjugate linear in the middle variable
and symmetric and bilinear in the outer variables satisfying that,
\begin{enumerate}[{\rm (a)}] \item $L(a,b) L(x,y) = L(x,y) L(a,b)
+ L(L(a,b)x,y)
 - L(x,L(b,a)y),$
where $L(a,b)$ is the operator on $E$ given by $L(a,b) x = \J
abx;$ \item $L(a,a)$ is an hermitian operator with non-negative
spectrum; \item $\|L(a,a)\| = \|a\|^2$.\end{enumerate}

For each $x$ in a JB$^*$-triple $E$, $Q(x)$ will stand for the
conjugate linear operator on $E$ defined by the assignment
$y\mapsto Q(x) y = \J xyx$.\smallskip

Every C*-algebra is a JB$^*$-triple via the triple product given
by $$2 \J xyz = x y^* z +z y^* x,$$ and every JB$^*$-algebra is a
JB$^*$-triple under the triple product \begin{equation}\label{e
triple product JB-alg} \J xyz = (x\circ y^*) \circ z + (z\circ
y^*)\circ x - (x\circ z)\circ y^*.
\end{equation}

The so-called \emph{Kaup-Banach-Stone} theorem for
JB$^{*}$-triples assures that a bounded linear surjection between
JB$^{*}$-triples is an isometry if and only if it is a triple
isomorphism (compare  \cite[Proposition 5.5]{Ka} or
\cite[Corollary 3.4]{BePe} or \cite[Theorem 2.2]{FerMarPe}). It
follows, among many other consequences, that when a
JB$^{*}$-algebra $J$ is a JB$^{*}$-triple for a suitable triple
product, then the latter coincides with the one defined in
$(\ref{e triple product JB-alg})$.\smallskip

A JBW$^*$-triple is a JB$^*$-triple which is also a dual Banach
space (with a unique isometric predual \cite{BarTi}). It is known
that the triple product of a JBW$^*$-triple is separately weak$^*$
continuous \cite{BarTi}. The second dual of a JB$^*$-triple $E$ is
a JBW$^*$-triple with a product extending the product of $E$
\cite{Di}.\smallskip

An element $e$ in a JB$^*$-triple $E$ is said to be a
\emph{tripotent} if $\J eee =e$. Each tripotent $e$ in $E$ gives
raise to the two following decomposition of $E$
$$E= E_{2} (e) \oplus E_{1} (e) \oplus E_0 (e),$$ where for
$i=0,1,2,$ $E_i (e)$ is the $\frac{i}{2}$ eigenspace of $L(e,e)$
(compare \cite[Theorem 25]{loos}). The natural projection of $E$
onto $E_i(e)$ will be denoted by $P_i(e)$. This decomposition is
termed the \emph{Peirce decomposition} of $E$ with respect to the
tripotent $e.$ The Peirce decomposition satisfies certain rules
known as \emph{Peirce arithmetic}: $$\J {E_{i}(e)}{E_{j}
(e)}{E_{k} (e)}\subseteq E_{i-j+k} (e),$$ if $i-j+k \in \{
0,1,2\}$ and is zero otherwise. In addition, $$\J {E_{2}
(e)}{E_{0}(e)}{E} = \J {E_{0} (e)}{E_{2}(e)}{E} =0.$$

The Peirce space $E_2 (e)$ is a JB$^*$-algebra with product
$x\circ_e y := \J xey$ and involution $x^{\sharp_e} := \J exe$.
\smallskip

A tripotent $e$ in $E$ is called \emph{complete} (resp., \emph{unitary})
if the equality $E_0(E)=0$ (resp., $E_2(e)=E$) holds.
When $E_2(e)=\CC e \neq \{0\},$ we say that $e$ is \emph{minimal}.\smallskip

For each element $x$ in a JB$^*$-triple $E$, we shall denote
$x^{[1]} := x$, $x^{[3]} := \J xxx$, and $x^{[2n+1]} := \J
xx{x^{[2n-1]}},$ $(n\in \NN)$. The symbol $E_x$ will stand for the
JB$^*$-subtriple generated by the element $x$. It is known that
$E_x$ is JB$^*$-triple isomorphic (and hence isometric) to $C_0
(\Omega)$ for some locally compact Hausdorff space $\Omega$
contained in $(0,\|x\|],$ such that $\Omega\cup \{0\}$ is compact,
where $C_0 (\Omega)$ denotes the Banach space of all
complex-valued continuous functions vanishing at $0.$ It is also
known that there exists a triple isomorphism $\Psi$ from $E_x$
onto $C_{0}(\Omega),$ satisfying $\Psi (x) (t) = t$ $(t\in \Omega)$ (cf.
\cite[Corollary 4.8]{Ka0}, \cite[Corollary 1.15]{Ka} and
\cite{FriRu}). The set $\overline{\Omega }=\hbox{Sp} (x)$ is
called the \emph{triple spectrum} of $x$. We should note that $C_0
(\hbox{Sp} (x)) = C(\hbox{Sp}(x))$, whenever $0\notin \hbox{Sp}
(x)$.\smallskip

Therefore, for each $x\in E$, there exists a unique element $y\in
E_x$ satisfying that $\J yyy =x$. The element $y,$ denoted by
$x^{[\frac13 ]}$, is termed the \emph{cubic root} of $x$. We can
inductively define, $x^{[\frac{1}{3^n}]} =
\left(x^{[\frac{1}{3^{n-1}}]}\right)^{[\frac 13]}$, $n\in \NN$.
The sequence $(x^{[\frac{1}{3^n}]})$ converges in the weak$^*$
topology of $E^{**}$ to a tripotent denoted by $r(x)$ and called
the \emph{range tripotent} of $x$. The tripotent $r(x)$ is the
smallest tripotent $e\in E^{**}$ satisfying that $x$ is positive
in the JBW$^*$-algebra $E^{**}_{2} (e)$ (compare \cite[Lemma
3.3]{EdRu}).\smallskip

A subspace $I$ of a JB$^*$-triple $E$ is a \emph{triple ideal} if
$\{E,E,I\}+\{E,I,E\} \subseteq I.$ By Proposition 1.3 in
\cite{BuChu}, $I$ is a triple ideal  if, and only if, the
inclusion $\{E,E,I\} \subseteq I$ holds. We shall say that $I$ is
an \emph{inner ideal} of $ E$ if $\{I, E, I\} \subseteq I.$ Given
an element $x$ in $E,$ let $E(x)$ denote the norm closed inner
ideal of $E$ generated by $x.$ It is known that $E(x)$ coincides
with the norm closure of the set $Q (x)(E).$ Moreover $E(x)$ is a
JB$^*$-subalgebra of $ E_2^{**}(r(x))$ and contains $x$ as a
positive element (compare \cite{BuZa}). Every triple ideal is, in
particular, an inner ideal. \smallskip

We recall that two elements $a,b$ in a JB$^*$-triple, $E,$ are
said to be \emph{orthogonal} (written $a\perp b$) if $L(a,b) =0$.
Lemma 1 in \cite{BurFerGarMarPe} shows that $a\perp b$ if and only
if one of the following nine statements holds:

\begin{equation}
\label{ref orthogo}\begin{array}{ccc}
  \J aab =0; & a \perp r(b); & r(a) \perp r(b); \\
  & & \\
  E^{**}_2(r(a)) \perp E^{**}_2(r(b));\ \ \ & r(a) \in E^{**}_0 (r(b));\ \ \  & a \in E^{**}_0 (r(b)); \\
  & & \\
  b \in E^{**}_0 (r(a)); & E_a \perp E_b & \J bba=0.
\end{array}
\end{equation}

The Jordan identity and the above reformulations assure that
\begin{equation}\label{orthog is subtriple} a\perp \J xyz,\hbox{
whenever } a\perp x,y,z.\end{equation}

An important class of JB$^*$-triples is given by the Cartan
factors. A JBW$^*$-triple $E$ is called factor if it contains no
proper weak$^*$ closed ideals. The Cartan factors are precisely
the JBW$^*$-triple factors containing a minimal tripotent
\cite{Ka97}. These can be classified in six different types (see
\cite{FriRu86} or \cite{Ka97}).\smallskip

A Cartan factor of type 1, denoted by $I_{n,m}$, is a
JB$^*$-triple of the form $L(H, H')$,  where $L(H,H')$ denotes the
space of bounded linear operators between two complex Hilbert
spaces $H$ and $H'$ of dimensions $n, m$ respectively, with the
triple product defined by $\{x, y, z\} = \frac{1}{2}(xy^*z +
zy^*x).$\smallskip

We recall that given a conjugation, $j$, on a complex Hilbert
space $H$, we can define the following linear involution $x\mapsto
x^t:=jx^*j$ on $L(H).$ A Cartan factor of type 2 (respectively,
type 3) denoted by $II_n$, (respectively, $III_n$) is the
subtriple of $L(H)$ formed by the $t$-skew-symmetric
(respectively, $t$-symmetric) operators, where $H$ is an $n$
dimensional complex Hilbert space. Moreover, $II_n$ and $III_n$
are, up to isomorphism, independent of the conjugation $j$ on
$H.$\smallskip

A Cartan factor of type 4, $IV_n$ (also called complex spin
factor),  is an $n$ dimensional  complex Hilbert space provided
with a conjugation $ x\mapsto \overline{x},$ where triple product
and the norm are  given by \begin{equation}\label{eq spin product}
\{x, y, z\} = (x|y)z + (z|y) x -(x|\overline{z})\overline{y},
\end{equation} and $ \|x\|^2 = (x|x) + \sqrt{(x|x)^2 -|
(x|\overline{x}) |^2},$ respectively.
\smallskip

The Cartan factor of type 6 is the 27-dimensional exceptional
JB$^*$-algebra $VI = H_3(\mathbb{O}^{\CC})$ of all symmetric 3 by
3 matrices with entries in the complex Octonions
$\mathbb{O}^{\CC}$, while the Cartan factor of type 5, $V =
M_{1,2}(\mathbb{O}^{\CC})$, is the subtriple of $
H_3(\mathbb{O}^{\CC})$ consisting of all $1\times 2$ matrices with
entries in $\mathbb{O}^{\CC}$.

\begin{rem}
\label{r min trips in spin}{\rm Let $E$ be a spin factor with
inner product $(.|.)$ and conjugation $x\mapsto \overline{x}$. It
is not hard to check (and part of the folklore of JB$^*$-triple
theory) that an element $w$ in $E$ is a minimal tripotent if and
only if $(w|\overline{w}) = 0$ and $(w|w)= \frac12$. For every
minimal tripotent $w$ in $E$ we have $E_2 (w) = \mathbb{C} w$,
$E_0 (w) = \mathbb{C} \overline{w}$ and $E_1 (w) = \{x\in E :
(x|w) = (x|\overline{w}) =0\}.$ Therefore, every minimal tripotent
$w_2\in E$ satisfying $w\perp w_2$ can be written in the form $w_2
= \lambda \overline{w}$ for some $\lambda\in \mathbb{C}$ with
$|\lambda|=1$. }
\end{rem}

\section{ Biorthogonality preservers }

Let $M$ be a subset of a JB$^*$-triple $E.$ We write
$M_{_E}^\perp$ for the \emph{(orthogonal) annihilator of $M$} defined by $$
M_{_E}^\perp:=\{ y \in E : y \perp x , \forall x \in M \}.$$ When
no confusion arise, we shall write $M^{\perp}$ instead of
$M^{\perp}_{_E}$.\smallskip

The next result summarises some basic properties of the
annihilator. The reader is referred to \cite[Lemma
3.2]{EdRu96} for a detailed proof.

\begin{lem}\label{l basic prop annihilator} Let $M$ a nonempty subset
 of a JB$^*$-triple $E$. \begin{enumerate}[{\rm $a)$}]
 \item $M^{\perp}$ is a norm closed inner ideal of $E$.
 \item $M\cap M^{\perp}= \{0\}.$
 \item $M\subseteq M^{\perp \perp}.$
 \item If $B\subseteq C$ then $C^{\perp}\subseteq B^{\perp}.$
 \item $M^{\perp}$ is weak$^*$ closed whenever $E$ is a
JBW$^*$-triple.$\hfill\Box$
 \end{enumerate}
\end{lem}

As illustration of the main identity (axiom $(a)$ in the definition of a JB$^*$-triple) 
we shall prove statement $a)$. For $a, a'$ in $M^{\perp}$, $b$ in $M$, and $c, d$ in $E$ 
we have $\J ca{\J d{a'}b} = \J {\J cad}{a'}b -
\J d{\J ac{a'}}b + \J d{a'}{\J cab}$, which shows that $\J ac{a'} \perp b$.\smallskip

Let $e$ be a tripotent in a JB$^*$-triple $E$. Clearly,
$\{e\}\subseteq E_2 (e)$. Therefore, it follows, by Peirce
arithmetic and Lemma \ref{l basic prop annihilator}, that $$E_2
(e)^{\perp} \subseteq \{e\}^\perp =E_0(e)\subseteq
E_2(e)^{\perp},$$ and hence
\begin{equation}\label{eq orth Peirce 2} E_2 (e)^{\perp} = \{e\}^\perp
=E_0(e).\end{equation} The next lemma describes the annihilator of
an element in an arbitrary JB$^*$-triple. Its proof follows
straightforwardly from the reformulations of orthogonality shown
in (\ref{ref orthogo}) (see also \cite[Lemma 1]{BurFerGarMarPe}).

\begin{lem} \label{l x.orth} Let $x$ be an element in a $JB^*$-triple $E$. Then
 $$\{x\}_{_E}^\perp =E^{**}_0(r(x))\cap E.$$ Moreover, when $E$ is
 a JBW$^*$-triple we have $$ \{x\}_{_{E}}^\perp =E_0(r(x)).$$ $\hfill \Box $
\end{lem}

\begin{proposition} \label{p orth.orth. trip}
Let  $e$ be a tripotent in a $JB^*$-triple $E$. Then
 $$ E_2(e) \oplus  E_1(e)\supseteq \{e\}_{_E}^{\perp \perp}= E_0(e)^{\perp} \supseteq E_2(e).$$
\end{proposition}

\begin{proof} It follows from (\ref{eq orth Peirce 2}) that $\{e\}^{ \perp \perp}=\{e\}^{ \perp \perp}_{_E}=(E_0(e))^\perp
\supseteq E_2(e).$ Now let us take $x \in (E_0(e))^\perp$. For
each $i\in \{0,1,2\}$ we write $x_i = P_i (e) (x)$, where $P_i
(e)$ denotes the Peirce $i$-projection with respect to $e$. Since
$x \in (E_0(e))^\perp,$ $x$ must be orthogonal to $x_0$ and so
$\{x_0,x_0,x\}=0.$ This equality, together with Peirce
arithmetic, show that $\{x_0,x_0,x_0\}+\{x_0,x_0,x_1\}=0,$ which
implies that $\|x_0\|^3 = \|\J {x_0}{x_0}{x_0} \|=0$.
\end{proof}

\begin{rem}
\label{r double annihilator}{\rm For a tripotent $e$ in a
JB$^*$-triple $E$, the equality $\{e\}_{_E}^{\perp \perp}=
E_0(e)^{\perp} = E_2(e)$, does not hold in general. Let $H_1$ and
$H_2$ be two infinite dimensional complex Hilbert spaces and let
$p$ be a minimal projection in $L(H_1)$. We define $E$ as the
orthogonal sum $p L(H_1) \oplus^{\infty} L(H_2)$. In this example
$\{p\}_{_E}^{\perp} = L(H_2)$ and $\{p\}_{_E}^{\perp\perp} =
pL(H_1)\neq \mathbb{C} p = E_2 (p).$\smallskip

 However, if $E$ is a Cartan factor and $e$ is a non-complete tripotent in $E,$ then the equality
 $\{e\}^{\perp \perp}=E_0(e)^\perp= E_2(e)$ always holds (compare Lemma 5.6 in \cite{Ka97}).}
\end{rem}

\begin{cor}\label{cor 4} Let  $x$ be an element in a $JB^*$-triple $E.$ Then

$$ E(x)\subseteq   E_2 ^{**}( r(x)) \cap E \subseteq \{ x\}_{_E}^{ \perp \perp
}.$$
\end{cor}

\begin{proof} Clearly, $E(x) = \overline{Q(x) (E)}\subseteq   E_2 ^{**}( r(x)) \cap
E$. Let us take $y$ in the intersection $E_2 ^{**}( r(x)) \cap E$.
In this case $y\in E_2 ^{**}( r(x)) \subseteq
\{x\}_{_{E^{**}}}^{\perp\perp}.$ Since $\{x\}_{_{E}}^{\perp}
\subset \{x\}_{_{E^{**}}}^{\perp}$, we conclude that $y\in
\{x\}_{_{E^{**}}}^{\perp\perp}\cap E \subseteq
\left(\{x\}_{_{E}}^{\perp}\right)_{E^{**}}^{\perp}\cap E =
\{x\}_{_{E}}^{\perp\perp}$.
\end{proof}

In the setting of C*-algebras the following sufficient conditions
to describe the second annihilator of a projection were
established in \cite[Lemma 3]{BurGarPe}.

\begin{lem}
\label{l bi annihilator closed projection}  Let $p$ be a
projection in a (non necessarily unital) C$^*$-algebra $A$. The
following assertions hold:
\begin{enumerate}[{\rm $ a)$}] \item $\{p\}_{_A}^{\perp} = (1-p) A (1-p)$,
where $1$ denotes the unit of $A^{**}$; \item
$\{p\}_{_A}^{\perp\perp} = p A p.$
\end{enumerate}$\hfill \Box$
\end{lem}

Let $x$ be an element in a JB$^*$-triple $E$. We say that $x$ is
 \emph{weakly compact} (respectively, \emph{compact}) if
the operator $Q(x):E\rightarrow E$ is weakly compact
(respectively, compact). A JB$^*$-triple is \emph{weakly compact}
(respectively, \emph{compact}) if every element in $E$ is weakly compact
(respectively, compact).\smallskip

Let $E$ be a JB$^*$-triple. If we denote by $K(E)$ the Banach
subspace of $E$ generated by its minimal tripotents, then $K(E)$
is a (norm closed) triple ideal of $E$ and it coincides with the
set of weakly compact elements of $E$ (see Proposition 4.7 in
\cite{BuChu}). For a Cartan factor $C$ we define the elementary
JB$^*$-triple of the corresponding type as $K(C).$ Consequently,
the elementary JB$^*$-triples, $K_i$ $(i = 1, ... , 6)$ are
defined as follows: $K_1 = K (H, H')$ (the compact operators
between two complex Hilbert spaces $H$ and $H'$); $K_i = C_i \cap
K(H)$ for $i = 2 , 3$, and $K_i = C_i$ for $i = 4, 5,
6$.\smallskip

It follows from \cite[Lemma 3.3 and Theorem 3.4]{BuChu} that a JB$^*$ triple, $E,$ is weakly compact if
and only if one of the following statement holds: \begin{enumerate}[$a)$] \item $K(E^{**})=K(E)$. \item $K(E)=E.$
\item $E$ is a $c_0$-sum of elementary JB$^*$-triples.
\end{enumerate}

Let $E$ be a JB$^*$-triple. A subset $S \subseteq E$ is said to be
\emph{orthogonal} if $0 \notin S$ and $x \perp y$ for every $x
\neq y$ in $S$. The minimal cardinal number $r$ satisfying
$card(S) \leq r$ for every orthogonal subset $S \subseteq E$ is
called the \emph{rank} of $E$ (and will be denoted by
$r(E)$).\smallskip

For every orthogonal family $(e_i)_{i\in I}$ of minimal tripotents
in a JBW$^*$-triple $E$ the w*-convergent sum $e: = \sum_{i} e_i$
is a tripotent and we call $(e_i)_{i\in I}$ a \emph{frame} in E if
$e$ is a maximal tripotent in $E$ (i.e., $e$ is a complete
tripotent and dim$(E_1 (e))\leq \hbox{dim} (E_1 (\widetilde{e}))$,
for every complete tripotent $\widetilde{e}$ in $E$). Every frame
is a maximal orthogonal family of minimal tripotents, the opposite
is not true in general (see \cite[\S 3]{BeLoPeRo} for more details).

\begin{proposition}
\label{p rank1} Let $e$ be a minimal tripotent in a $JB^*$-triple
$E$. Then $\{e\}_{_E}^{\perp \perp}$ is a rank-one norm closed
inner ideal of $E$.
\end{proposition}

\begin{proof} Let $F$ denote $\{e\}_{_E}^{\perp\perp}$. Since $e$
is a minimal tripotent (i.e. $E_{2} (e) = \mathbb{C} e$), the set of
states on $E_{2} (e)$, $\{ \varphi\in E^* : \varphi (e) =1 = \|\varphi\|\}$,
reduces to one point $\varphi_0$ in $E^*$. Proposition 2.4 and Corollary 2.5 in
\cite{BuFerMarPe} imply that the norm of $E$ restricted to $E_1 (e)$ is equivalent
to a Hilbertian norm. More precisely, in the terminology of \cite{BuFerMarPe}, the norm
$\|.\|_{e}$ coincides with the Hilbertian norm $\|.\|_{\varphi_0}$ and
is equivalent to the norm of $E_1 (e)$.\smallskip

Proposition \ref{p orth.orth. trip} guarantees that $F$
is a norm closed subspace of $E_2 (e) \oplus E_1 (e) = \mathbb{C}
e \oplus E_1 (e)$, and hence $F$ is isomorphic to a Hilbert
space.\smallskip

We deduce, by Proposition 4.5.(iii)  in \cite{BuChu} (and its
proof), that $F$ is a finite orthogonal sum of Cartan factors,
$C_1,\ldots, C_m$, which are finite-dimensional, or infinite
dimensional spin factors, or of the form $L(H,H')$ for suitable
complex Hilbert spaces $H$ and $H'$ with dim$(H') < \infty$. Since
$F$ is an inner ideal of $E$ (and hence a JB$^*$-subtriple of $E$)
and $e$ is a minimal tripotent in $E$, we can easily check that $e$ is a
minimal tripotent in $F=\bigoplus^{\ell_{\infty}}_{j=1,\ldots,m} C_j$.
If we write $e= e_1+\ldots+e_m$, where each $e_j$ is a tripotent in $C_j$
and $e_j\perp e_k$ whenever $j\neq k$, since $\mathbb{C} e_1\oplus\ldots
\oplus \mathbb{C} e_1\subseteq F_{2} (e) = \mathbb{C} e$, we deduce that
there exists a unique $j_0\in \{1,\ldots,m\}$ satisfying $e_j =0$ for
all $j\neq j_0$ and $e=e_{j_0}\in C_{j_0}$.\smallskip

For each $j\neq j_0,$ we have $C_{j}\subseteq \{e\}_{E}^{\perp},$ and hence
$$\bigoplus^{\ell_{\infty}}_{j=1,\ldots,m} C_j = F = \{e\}^{\perp\perp}
\subseteq C_{j}^{\perp}.$$ This implies that $C_j\perp C_j$ (or equivalently, $C_j=0$)
for every $j\neq j_0$. We consequently have $F =
\{e\}_{_E}^{\perp\perp} = C_{j_0}$.\smallskip

Finally, if $r(F)\geq  2$, then we deduce, via Proposition 5.8 in
\cite{Ka97}, that there exists minimal tripotents
$e_2,\ldots,e_{r}$ in $ F$ satisfying that $e,e_2,\ldots,e_r$ is a
frame in $F$. For each $i\in \{2,\ldots,r\}$, $e_i $ is orthogonal
to $e$ and lies in $F=\{e\}_{_E}^{\perp\perp}$, which is
imposible.
\end{proof}

Let  $T: E\rightarrow F$ be a linear map between two
JB$^*$-triples. We shall say that $T$ is \emph{orthogonality
preserving} if $T(x) \perp T(y)$ whenever $x \perp y$. The mapping
$T$ is said to be \emph{biorthogonality preserving} whenever the
equivalence $$ x \perp y \Leftrightarrow T(x) \perp T(y)$$ holds
for all $x,y$ in $E$.\smallskip

It can be easily seen that every biorthogonality preserving
linear mapping $T: E\to F,$ between JB$^*$-triples is
injective. Indeed, for each $x\in E$, the condition $T(x)=0$
implies that $T(x)\perp T(x)$, and hence $x\perp x$, which gives
$x=0$.\smallskip

Orthogonality preserving bounded linear maps from a JB$^*$-algebra
to a JB$^*$-triple were completely described in
\cite{BurFerGarPe}.\smallskip

Before stating the result, let us recall some basic definitions.
Two elements $a$ and $b$ in a JB$^*$-algebra $J$ are said to
\emph{operator commute} in $J$ if the multiplication operators
$M_a$ and $M_b$ commute, where $M_a$ is defined by $M_a (x) :=
a\circ x$. That is, $a$ and $b$ operator commute if and only if
$(a\circ x) \circ b = a\circ (x\circ b)$ for all $x$ in $J$.
Self-adjoint elements $a$ and $b$ in $J$ generate a
JB$^*$-subalgebra that can be realised as a JC$^*$-subalgebra of
some $B(H)$, \cite{Wri77}, and, in this realisation, $a$ and $b$
commute in the usual sense whenever they operator commute in $J$
\cite[Proposition 1]{Top}.  Similarly, two self-adjoint elements
$a$ and $b$ in $J$ operator commute if and only if $a^2 \circ b=
\J aab =\J aba$ (i.e., $a^2 \circ b = 2 (a\circ b)\circ a - a^2
\circ b$). If $b\in J$ we use $\{b\}'$ to denote the set of
elements in $J$ that operator commute with $b$. We shall write
$Z(J):= J'$ for the center of $J$ (This agrees with the usual
notation in von Neumann algebras).\smallskip

\begin{thm}\label{t BurFerGarPe}\cite[Theorem 4.1 ]{BurFerGarPe}  Let $T: J\to E$ be a
bounded linear mapping from a JB$^*$-algebra to a JB$^*$-triple.
For $h= T^{**} (1)$ and $r= r(h)$ the following assertions are
equivalent.
\begin{enumerate}[{\rm $a)$}]
 \item $T$ is orthogonality
preserving.

\item There exists a unique Jordan $^*$-homomorphism $S: J \to
E_2^{**} (r) $ satisfying that $S^{**} (1) = r$, $S(J)$ and $h$
operator commute and $T(z) = h\circ_r S(z) $ for all $z\in J$.

\item $T$ preserves zero triple products, that is, $\J
{T(x)}{T(y)}{T(z)} = 0$ whenever $\J xyz =0$.$\hfill\Box$
\end{enumerate}
\end{thm}

The above characterisation proves that the bitranspose of an
orthogonality preserving bounded linear mapping from a
JB$^*$-algebra onto a JB$^*$-triple also is orthogonality
preserving.\smallskip

The following theorem was essentially proved in \cite{BurFerGarPe}.
We include here an sketch of the proof for completeness reasons.

\begin{thm}\label{charbbiops} Let $T:J \rightarrow E$ be a surjective linear
operator from a  JBW$^*$-algebra onto a JBW$^*$-triple and let $h$
denote $T(1)$. Then $T$ is biorthogonality preserving if, and only
if, $r(h)$ is a unitary tripotent in $E,$ $h$ is an invertible element in the
JB$^*$-algebra $E=E_2(r(h)),$ and there exists a Jordan
*-isomorphism $S:J \rightarrow E=E_2(r(h))$ such that
$S(J)\subseteq \{ h\}'$ and $T=h\circ_{r(h)}S.$ Further, when $J$
is a factor {\rm(}i.e. $Z(J) = \mathbb{C} 1${\rm)} then $T$ is a scalar
multiple of a triple isomorphism.
 \end{thm}

\begin{proof} The sufficient implication is clear. We shall prove the necessary condition.
To this end let $T:J \rightarrow E$ be a surjective linear
operator from a JBW$^*$-algebra onto a JBW$^*$-triple and let $h=T(1)\in E$. We have already
commented that every biorthogonality preserving linear mapping between JB$^*$-triples is
injective. Therefore $T$ is a linear bijection.\smallskip

From Corollary 4.1 $(b)$ and its proof in \cite{BurFerGarPe}, we deduce that $$T(J_{sa}) \subseteq E_2
(r(h))_{sa}, \hbox{ and hence } E=T(J) \subseteq E_2 (r(h))\subseteq E.$$ This implies that $E=E_2 (r(h))$,
which assures that $r(h)$ is a unitary tripotent in $E$. Since the range tripotent of $h$, $r(h)$, is the
unit of $E_2 (r(h))$ and $h$ is a positive element in the JBW$^*$-algebra $E_2 (r(h)),$ we can easily check that
$h$ is invertible in $E_2 (r(h))$. Furthermore, $h^{^{\frac12}}$ is invertible in $E_2 (r(h))$ with inverse
$h^{^{-\frac12}}$.\smallskip

The proof of \cite[Theorem 4.1]{BurFerGarPe}
can be literally applied here to show the existence of a Jordan *-homomorphism $S:J \rightarrow E=E_2(r(h))$
such that $S(J)\subseteq \{ h\}'$ and $T=h\circ_{r(h)} S.$ Since, for each $x\in J$, $h$ and $S(x)$ operator
commute and $h^{\frac12}$ lies in the JB$^*$-subalgebra of $E_2(r(h))$ generated by $h$, we can easily check
that $S(x)$ and $h^{\frac12}$ operator commute. Thus, $$T=h\circ_{r(h)} S = U_{h^{^{\frac12}}} S,$$ where
$U_{h^{^{\frac12}}}: E_2(r(h)) \to E_2(r(h))$ is the linear mapping defined by $$U_{h^{^{\frac12}}} (x)
= 2 ({h^{^{\frac12}}}\circ_{r(h)} x) \circ_{r(h)} {h^{^{\frac12}}} - ({h^{^{\frac12}}}\circ_{r(h)} {h^{^{\frac12}}}) \circ_{r(h)} x.$$
It is well known that ${h^{^{\frac12}}}$ is invertible if and only if $U_{h^{^{\frac12}}}$ is an invertible operator and,
in this case, $U_{h^{^{\frac12}}}^{-1} = U_{h^{^{-\frac12}}}$ (compare \cite[Lemma 3.2.10]{HoS}). Therefore, $S = U_{h^{-\frac12}} T.$
It follows from the bijectivity of $T$ that $S$ is a Jordan *-isomorphism.\smallskip

Finally, when $Z(J) = \mathbb{C} 1$, the center of $E_2(r(h))$ also reduces to $\mathbb{C} r(h)$ and,
since $h$ is an invertible element
in the center of $E_2(r(h))$, we deduce that $T$ is a scalar multiple of a triple isomorphism.
\end{proof}

\begin{proposition}
\label{p direct sums} Let $E_1$, $E_2$ and $F$ be three
JB$^*$-triples (respectively, JBW$^*$-triples). Let us suppose that $T : E_1\oplus^{\infty} E_2
\to F$ is a biorthogonality preserving linear surjection. Then
$T(E_1)$ and $T(E_2)$ are norm closed (respectively, weak$^*$ closed)
inner ideals of $F$, $B=T(A_1)\oplus^{\infty} T(A_2),$ and for $j= 1,2$ $T|_{A_j} : A_j
\to T(A_j)$ is a biorthogonality preserving linear surjection.
\end{proposition}

\begin{proof}
Let us fix $j\in \{1,2\}$. Since $E_{j} = E_{j}^{\perp\perp}$ and
$T$ is a biorthogonality preserving linear surjection, we deduce
that $T(E_j) = T(E_j^{\perp\perp}) = T(E_j)^{\perp\perp}$. Lemma
\ref{l basic prop annihilator} guarantees that $T(E_j)$ is a norm
closed inner ideal of $F$ (respectively, a weak$^*$ closed inner ideal of
$F$ whenever $E_1$, $E_2$ and $F$ are JBW$^*$-triples). The
rest of the proof follows from Lemma \ref{l basic prop annihilator} and the fact
that $F$ coincides with the orthogonal sum of $T(E_1)$ and $T(E_2)$.
\end{proof}

\section{Biorthogonality preservers between weakly compact JB$^*$-triples}

The following theorem generalises \cite[Theorem 5]{BurGarPe} by
proving that biorthogonality preserving linear surjections between
JB$^*$-triples send minimal tripotents to scalar multiples of
minimal tripotents.

\begin{thm}\label{thm minimal.to.minimal} Let $T: E\rightarrow F$ be a
biorthogonality preserving linear surjection between two
JB$^*$-triples and let $e$ be a minimal tripotent in $E$. Then
$\|T(e)\|^{-1} \ T (e) = f_e$ is a minimal tripotent in $F$.
Further, $T$ satisfies that $T (E_2 (e)) = F_2 (f_e) $ and $T(E_0
(e)) = F_0 (f_e)$.
\end{thm}

\begin{proof} Since $T$ is biorthogonality preserving surjection, the equality
$$T(S_{_E}^\perp) = T(S)_{_F}^\perp$$ holds for every subset $S$ of
$E.$ Lemma \ref{l basic prop annihilator} assures that for each
minimal tripotent $e$ in $E$, $\{T(e)\}_{_F}^{ \perp \perp } =
T(\{e\}_{_E}^{ \perp \perp })$ is a norm closed inner ideal in
$F$. By Proposition \ref{p rank1}, $\{e\}_{_E}^{ \perp \perp }$
is a rank-one JB$^*$-triple, and hence $\{T(e)\}_{_F}^{
\perp \perp }$ cannot contain two non-zero orthogonal elements.
Thus, $\{T(e)\}_{_F}^{\perp \perp }$ is a rank-one JB$^*$-triple.\smallskip

The arguments given in the proof of Proposition \ref{p rank1}
above (see also Proposition 4.5.(iii) in \cite{BuChu} and its
proof or \cite[\S 3]{BeLoPeRo}) show that the inner ideal $\{T(e)\}_{_F}^{\perp \perp }$ is
a rank-one Cartan factor, and hence a type 1 Cartan factor of the
form $L(H,\mathbb{C})$, where $H$ is a complex Hilbert space or a
type 2 Cartan factor $II_{3}$ (it is known that $II_{3}$ is
JB$^*$-triple isomorphic to a 3-dimensional complex Hilbert
space). This implies that $\|T(e)\|^{-1} \ T (e) = f_e$ is a
minimal tripotent in $F$ and $ T(e)= \lambda_{e}\ f_e,$ for a
suitable $\lambda_{e} \in \CC \setminus \{0\}.$\smallskip

The equality $T (E_2 (e)) = F_2 (f_e) $ has been proved.
Concerning the Peirce zero subspace we have $$T(E_0 (e)) = T(E_2
(e)_{_{E}}^{\perp}) = T\left(E_2 (e)\right)_{_F}^{\perp} = F_2
(f_e)_{_F}^{\perp} = F_0 (f_e).$$
\end{proof}\smallskip

Let $H$ and $H'$ be complex Hilbert spaces. Given $k\in H'$ and
$h\in H$, we define $k\otimes h$ in $L(H,H')$ by the assignment
$k\otimes h (\xi) := (\xi | h) \ k$. In this case, every minimal
tripotent in $L(H,H')$ can be written in the form $k\otimes h$,
where  $h$ and $k$ are norm-one elements in $H$ and $H'$,
respectively. It can be easily seen that two minimal tripotents
$k_1 \otimes h_1$ and $ k_2\otimes h_2$ are orthogonal if, and
only if, $h_1 \perp h_2$ and $k_1\perp k_2.$\smallskip

\begin{thm}\label{biopcartan}
Let $T:E \rightarrow F$ be a biorthogonality preserving linear
surjection between two JB$^*$-triples, where $E$ is a type
$I_{n,m}$ Cartan factor, with $n,m\geq 2 $. Then there exists  a
positive real number $\lambda$ such that $\|T(e)\|=\lambda,$ for
every minimal tripotent $e$ in $E$.
\end{thm}

\begin{proof} Let $H,H'$ be complex Hilbert spaces such that $E=L(H,H')$.
Let $e_1:= k_1 \otimes h_1$ and $e_2 := k_2\otimes h_2$ be two
minimal tripotents in $E$. We write $H_1=$Span $(\{h_1,h_2\})$ and
$H'_1=$Span$( \{k_1,k_2\}).$ The tripotents $k_1 \otimes h_1$ and
$ k_2\otimes h_2$ can be identified with elements in
$L(H_1,H'_1)$. By Theorem \ref{thm minimal.to.minimal},
$T(e_1)=\alpha_1 \ f_1$ and $T(e_2)=\alpha_2 \ f_2$, where $f_1$ and
$f_2$ are two minimal tripotents in $F$.\smallskip

If dim$(H_1)=\hbox{dim}(H'_1)=2,$ then the norm closed inner
ideal, $E_{e_1,e_2}$, of $E$ generated by $e_1$ and $e_2$
identifies with $L(H_1,H'_1),$ which is JB$^*$-isomorphic to
$M_2(\CC)$ and coincides with the inner ideal generated by two orthogonal
minimal tripotents $g_1=\left(
                          \begin{array}{cc}
                            1 & 0 \\
                            0 & 0 \\
                          \end{array}
                        \right)
$ and $g_2= \left(
              \begin{array}{cc}
                0 & 0 \\
                0 & 1 \\
              \end{array}
            \right)
$, where $g_1+g_2$ is the unit element in $E_{e_1,e_2}\cong M_2 (\mathbb{C})$.\smallskip

By Theorem \ref{thm minimal.to.minimal},
$w_1:=\frac{1}{\|T(g_1)\|} T(g_1)$ and $w_2:=\frac{1}{\|T(g_2)\|}
T(g_2)$ are two orthogonal minimal tripotents in $F$. The element
$w= w_1 +w_2$ is a rank-2 tripotent in $F$ and coincides with the
range tripotent of the element $h=T(g_1+g_2)= \|T(g_1)\| w_1 + \|T(g_2)\| w_2$. By Theorem \ref{t
BurFerGarPe} (see also \cite[Corollary 4.1 $(b)$]{BurFerGarPe}),
${T(E_{e_1,e_2})} \subseteq F_{2} (w)$. It is not hard to see that
$h$ is invertible in $F_{2} (w)$ with inverse $h^{-1}=\frac{1}{\|T(g_1)\|}w_1+\frac{1}{\|T(g_2)\|}w_2$.\smallskip

The inner ideal $E_{e_1,e_2}$ is finite dimensional, $T(E_{e_1,e_2})$ is norm closed and
$T|_{E_{e_1,e_2}}: E_{e_1,e_2} \to F$ is a continuous
biorthogonality preserving linear operator. Theorem \ref{t BurFerGarPe} guarantees
the existence of a Jordan $^*$-homomorphism $S: E_{e_1,e_2}\cong M_2 (\mathbb{C}) \to
F_2 (w) $ satisfying that $S (g_1+g_2) = w$, $S(E_{e_1,e_2})$ and $h$
operator commute and \begin{equation}\label{eq new} T(z) = h\circ_w S(z), \hbox{ for all } z\in E_{e_1,e_2}.
\end{equation} It follows
from the operator commutativity of $h^{-1}$ and $S(E_{e_1,e_2})$ that $S(z) = h^{-1}\circ_w T(z)$
for all $z\in E_{e_1,e_2}$. The injectivity of $T$ implies that $S$ is a Jordan *-monomorphism.\smallskip

Lemma 2.7 in \cite{FerMarPe2}
shows that $F_2 (w) = F_2 (w_1 +w_2)$ coincides with
$\mathbb{C}\oplus^{\ell_{\infty}}\mathbb{C}$ or with a spin
factor. Since $4= \dim ({T(E_{e_1,e_2})}) \leq \dim(F_{2} (w))$,
we deduce that $F_2 (w)$ is a spin factor with inner product
$(.|.)$ and conjugation $x\mapsto \overline{x}$. From Remark
\ref{r min trips in spin}, we may assume, without loss of
generality, that $(w_1|w_1) = \frac12,$ $(w_1|\overline{w}_1)=0$,
and $w_2 = \overline{w}_1$.\smallskip

Now, we take $g_3=\left(
\begin{array}{cc}
0 & 1 \\
0 & 0 \\
\end{array}
\right) $ and $g_4= \left(
              \begin{array}{cc}
                0 & 0 \\
                1 & 0 \\
              \end{array}
            \right)$ in $E_{e_1,e_2}$. The elements $w_3:= S(g_3)$ and $w_4:= S(g_4)$
are two orthogonal minimal tripotents in $F_2 (w)$ with $\J
{w_i}{w_i}{w_j} = \frac12 w_j$ for every $(i,j),(j,i) \in
\{1,2\}\times \{3,4\}$. Applying again Remark \ref{r min trips in
spin}, we may assume that $(w_3|w_3) = \frac12,$
$(w_3|\overline{w}_3)=0$, $w_4 = \overline{w}_3$, and $(w_3|w_1) =
(w_3|w_2)=0$. Applying the definition of the triple product in a
spin factor given in $(\ref{eq spin product})$ we can check that
$(w_1,w_3,w_2=\overline{w}_1,w_4=\overline{w}_3)$ are four minimal
tripotents in $F_{2} (w)$ with $w_1\perp w_2$, $w_3\perp w_4$, $\J
{w_i}{w_i}{w_j} = \frac12 w_j$ for every $(i,j),(j,i) \in
\{1,2\}\times \{3,4\}$, $\J {w_1}{w_3}{w_2} = -\frac12 w_4$, $\J
{w_3}{w_2}{-w_4} = \frac12 w_1$, $\J
{w_2}{-w_4}{w_1} = \frac12 w_3$, and $\J
{-w_4}{w_1}{w_3} = \frac12 w_2$.
Therefore, denoting by $M$ the JB$^*$-subtriple of $F_2 (w)$ generated by
$w_1,w_3,w_2,$ and $w_4$, we have shown that $M$ is JB$^*$-triple isomorphic
to $M_2 (\mathbb{C})$.\smallskip

Combining $(\ref{eq new})$ and $(\ref{eq spin product})$ we get
$$T(g_3) = h\circ_{w} S(g_3) = \J hw{w_3} = \frac{\|T(g_1)\|+\|T(g_2)\|}{2} w_3$$ and
$$T(g_4) = h\circ_{w} S(g_4) = \J hw{w_4} = \frac{\|T(g_1)\|+\|T(g_2)\|}{2} w_4.$$\smallskip
Since $T(g_1) = \|T(g_1)\|\ w_1$, $T(g_2) = \|T(g_2)\|\ w_2$, and $E_{e_1,e_2}$ is linearly
generated by $g_1,g_2,g_3$ and $g_4,$ we deduce that $T(E_{e_1,e_2}) \subseteq M$ with
$4= \dim(T(E_{e_1,e_2})) \leq  \dim(M) = 4$. Thus, $T(E_{e_1,e_2}) = M$ is a JB$^*$-subtriple of $F$.\smallskip

The mapping $T|_{E_{e_1,e_2}}: E_{e_1,e_2}\cong M_2(\mathbb{C}) \to {T(E_{e_1,e_2})} $ is
a continuous biorthogonality preserving linear bijection. Theorem \ref{charbbiops}
implies that $T|_{E_{e_1,e_2}}$ is a (non-zero) scalar multiple of a triple isomorphism,
and hence $\|T(e_1)\|=\|T(e_2)\|.$\smallskip

If $dim(H'_1)=1,$ then $L(H_1,H'_1)$ is a rank-one JB$^*$-triple.
Since $n,m\geq 2$, we can find a minimal tripotent $e$ in $ E$
such that the norm closed inner ideals of $E$ generated by $\{e,
e_1\}$ and $\{e,e_2\}$ both coincide with $M_2(\CC).$
The arguments given in the above paragraph show that
$\|T(e_1)\|=\|T(e)\|=\|T(e_2)\|.$\smallskip

Finally, the case $dim(H_1)=1$ follows from
the same arguments given in the above paragraph.
\end{proof}

\begin{rem}
\label{r Cauchy series}{\rm Given a sequence $(\mu_n)\subset c_0$
and a bounded sequence $(x_n)$ in a Banach space $X$, the series
$\sum_{k} \mu_k x_k$ need not be, in general, convergent in $X$.
However, when $(x_n)$ is a bounded sequence of mutually orthogonal
elements in a JB$^*$-triple, $E$, the equality
$$\left\| \sum_{k=1}^{n} \mu_k x_k - \sum_{k=1}^{m} \mu_k x_k \right\| = \max \left\{ |\mu_{n+1}|,\ldots , |\mu_{m}|\right\} \sup\{\|x_n\|\},$$ holds for every $n<m$ in $\mathbb{N}$. It follows that $(\sum_{k=1}^{n} \mu_k x_k)$ is a Cauchy series and hence convergent in $E$.}
\end{rem}

The following three results generalise \cite[Lemmas 8, 9 and
Proposition 10]{BurGarPe} to the setting of JB$^*$-triples.

\begin{lem}\label{l c0 sums} Let $T: E\rightarrow F$ be a
biorthogonality preserving linear surjection between two
JB$^*$-triples and let $(e_n)_{n}$ be a sequence of mutually
orthogonal minimal tripotents in $E$. Then there exist two
positive constants $ m\leq M $ satisfying $m \leq \|T(e_n)\| \leq
M$, for all $n\in \mathbb{N}$.
\end{lem}

\begin{proof} We deduce from Theorem \ref{thm minimal.to.minimal} that,
for each natural $n$, there exist a minimal tripotent $f_n$ and a
scalar $\lambda_n\in \mathbb{C}\setminus \{0\}$ such that $T(e_n)
= \lambda_n f_n$, where $\|T(e_n)\| =  \lambda_n$. Note that $T$
being biorthogonality preserving implies $(f_n)$ is a sequence of
mutually orthogonal minimal tripotents in $F$.\smallskip

Let $(\mu_n)$ be any sequence in $c_0$. Since the elements $e_n$'s
are  mutually orthogonal the series $\sum_{k\geq 1}  \mu_k e_k$
converges to an element in $E$ (compare Remark \ref{r Cauchy
series}). For each natural $n$, $\sum_{k\geq 1}^{\infty}  \mu_k
e_k $ decomposes as the orthogonal sum of $\mu_n e_n$ and
$\sum_{k\neq n}^{\infty}  \mu_k e_k $, therefore
$$T\left(\sum_{k\geq 1}^{\infty}  \mu_k e_k \right) = \mu_n
\lambda_n f_n  + T\left(\sum_{k\neq n}^{\infty}  \mu_k
e_k\right),$$ with $\mu_n \lambda_n f_n  \perp T\left(\sum_{k\neq
n}^{\infty}  \mu_k e_k\right)$, which in particular implies
$$\left\|T\left(\sum_{k\geq 1}^{\infty}  \mu_k e_k \right)
\right\|= \max\left\{|\mu_n| | \lambda_n|, \left\|
T\left(\sum_{k\neq n}^{\infty}  \mu_k
e_k\right)\right\|\right\}\geq |\mu_n| | \lambda_n|. $$ This
establishes that, for each $(\mu_n)$ in $c_0$, $(\mu_n
\lambda_n)$ is a bounded sequence, which in particular implies
that $(\lambda_n)$ is bounded.\smallskip

Finally, since $T$ is a biorthogonality preserving linear
surjection and \linebreak $T^{-1} (f_n) = \lambda_n^{-1} e_n$,  we
can similarly show that $(\lambda_n^{-1})$ is also bounded.
\end{proof}

\begin{lem}\label{l lim series to zero}
Let $T: E\rightarrow F$ be a biorthogonality preserving linear
surjection between two JB$^*$-triples, $(\mu_n)$ a sequence in
$c_0$ and let $(e_n)_{n}$ be a sequence of mutually orthogonal
minimal tripotents in $E$. Then the sequence $\left( T \left(
\sum_{k\geq n}^{\infty}  \mu_k e_k \right)\right)_{n}$ is well
defined and converges in norm to zero.
\end{lem}

\begin{proof} By Theorem \ref{thm minimal.to.minimal} and Lemma \ref{l c0 sums}
it follows that $(T(e_n))$ is a bounded sequence of mutually
orthogonal elements  in $F$. Let $M$ denote a bound of the above
sequence. For each natural $n$, Remark \ref{r Cauchy series}
assures that the series $\sum_{k\geq n}^{\infty}  \mu_k e_k$
converges.\smallskip

Let us define $y_n := T \left( \sum_{k\geq n}^{\infty}  \mu_k e_k
\right)$. We claim that $(y_n)$ is a Cauchy sequence in $F$.
Indeed, given $n< m$ in  $\mathbb{N}$, we have \begin{equation}
\label{eq lemma lim series to zero} \| y_n -y_m\| = \left\| T
\left( \sum_{k\geq n}^{m-1}  \mu_k e_k \right) \right\| = \left\|
\sum_{k\geq n}^{m-1}  \mu_k T(e_k) \right\|
\end{equation} $$ \leq M\ \max \{ |\mu_{n}|,\ldots , |\mu_{m-1}|\},$$
where in the last inequality we used the fact that $(T(e_n))$ is a
sequence  of mutually orthogonal elements. Consequently, $(y_n)$
converges in norm to some element $y_0$ in $F$. Let $z_0$ denote
$T^{-1} (y_0)$.\smallskip

Let us fix a natural $m$. By hypothesis, for each $n>m$, $e_m$ is
orthogonal to $\sum_{k\geq n}^{\infty}  \mu_k e_k$. This implies
that $T(e_m) \perp y_n$, for every $n>m$, which, in particular,
implies $\J {T(e_m)}{T(e_m)}{y_n} =0,$ for every $n>m$. Taking
limit when $n$ tends to $\infty$ we have $\J {T(e_m)}{T(e_m)}{y_0}
=0.$ This shows that $y_0=T(z_0)$ is orthogonal to $T(e_m)$, and
hence $e_m\perp z_0$. Since $m$ was arbitrarily chosen, we deduce
that $z_0$ is orthogonal to $\sum_{k\geq n}^{\infty}  \mu_k e_k,$
for every natural $n$. Therefore, $(y_n)\subset \{ y_0\}^{\perp},$
and hence $y_0 $ belongs to the norm closure of $\{
y_0\}^{\perp},$ which implies $y_0 =0$.
\end{proof}

\begin{proposition}
\label{t weakly compact} Let $T: E\rightarrow F$ be a
biorthogonality preserving linear surjection between two
JB$^*$-triples, where $E$ is weakly compact. Then $T$ is
continuous if and only if the set $\mathcal{T}:=\left\{ \|T(e)\| :
e \hbox{ minimal tripotent in } E\right\} $ is bounded. Moreover,
if that is the case, then $\|T\|=\sup (\mathcal{T}).$
\end{proposition}

\begin{proof}
The necessity being obvious. Suppose that $$M= \sup \left\{
\|T(e)\| : e \hbox{ minimal tripotent in } E\right\} <\infty.$$
Since $E$ is weakly compact, each nonzero element $x$ of $E$ can
be written as  a norm convergent (possibly finite) sum $x =
\sum_{n} \lambda_n u_n$, where $u_n$ are mutually orthogonal
minimal tripotents of $E$, and $\|x\| = \sup \{|\lambda_n| :
n\}$(compare Remark 4.6 in \cite{BuChu}). If the series $x=
\sum_{n} \lambda_n u_n$ is finite then $$\|T(x) \| = \left\|
\sum_{n=1}^{m} \lambda_n T(u_n) \right\| \stackrel{(*)}{=} \max
\left\{  \left\| \lambda_n T(u_n) \right\| : n= 1,\ldots,m\right\}
\leq M \|x\|,$$ where at $(*)$ we apply the fact that $(T(u_n))$
is a finite set of mutually orthogonal tripotents in $F$. When the
series $x=\sum_{n} \lambda_n u_n$ is infinite we may assume that
$(\lambda_n)\in c_0$.\smallskip

It follows from Lemma \ref{l lim series to zero} that  the
sequence $\left( T \left( \sum_{k\geq n}^{\infty}  \lambda_k u_k
\right)\right)_{n}$ is well defined and converges in norm to zero.
We can find a natural $m$ such that $\left\|T \left( \sum_{k\geq
m}^{\infty}  \lambda_k u_k \right)\right\| < M \|x\|$. Since the
elements $\lambda_1 u_1,\ldots,\lambda_{m-1} u_{m-1},$ $
\sum_{k\geq m}^{\infty}  \lambda_k u_k$ are mutually orthogonal,
we have $$\|T(x) \| = \max \big\{ \|T(\lambda_1
u_1)\|,\ldots,\|T(\lambda_{m-1} u_{m-1})\|, \left\|T(\sum_{k\geq
m}^{\infty}  \lambda_k u_k)\right\| \big\} \leq M \|x\|.$$
\end{proof} \smallskip

Let $E$ be an elementary JB$^*$-triple of type 1 (that is,  an
elementary JB$^*$-triple such that $E^{**} $ is a type $1$ Cartan
factor), and let $T:E \rightarrow F$ a biorthogonality preserving
linear surjection from $E$ to another JB$^*$-triple. Then by
Theorem \ref{biopcartan} and Proposition \ref{t weakly compact},
$T$ is continuous. Further, we claim that $T$ is a scalar multiple
of a triple isomorphism. Indeed, let us see that
$S=\frac{1}{\lambda}T$ is a triple isomorphism, where
$\lambda=\|T(e)\|=\|T\|$ for some (and hence any) minimal
tripotent $e$ in $E$ (compare Theorem \ref{biopcartan}). Let $x$
be an element of $E, $ then $x=\sum_n \lambda_n e_n$ for suitable
$(\lambda_n) \in c_0$ and a family of mutually orthogonal minimal
tripotents $(e_n)$ in $E$ \cite[Remark 4.6]{BuChu}. Then by
observing that $T$ is continuous we have

$$\|S(x)\|=\frac{1}{\lambda}\|T(x)\|=\frac{1}{\lambda}\left\|T(\sum_n \lambda_n e_n)\right\|$$
$$=\frac{1}{\lambda}\left\|\sum_n \lambda_n T(e_n)\right\|=\frac{1}{\lambda}\sup_n
\{ |\lambda_n| \ \|T(e_n)\| \}$$
$$=\frac{1}{\lambda}\sup_n \{ |\lambda_n| \ \lambda \}=\sup_n \{ |\lambda_n| \}=\|x\|.$$
This proves that $S$ is a
surjective linear isometry between JB$^*$-triples, and hence a
triple isomorphism (see \cite[Proposition 5.5]{Ka} or
\cite[Corollary 3.4]{BePe} or \cite[Theorem 2.2]{FerMarPe}). We
have thus proved the following result:

\begin{cor}\label{scalartripletypeI} Let $T:E \rightarrow F$ a
biorthogonality preserving linear surjection from a type 1
elementary JB$^*$-triple of rank bigger than one, onto another
JB$^*$-triple. Then $T$ is a scalar multiple of a triple
isomorphism. $\hfill \Box$
\end{cor}\smallskip

Let $p$ and $q$ be two minimal projections in a C$^*$-algebra $A$
with $q\neq p$. It is known that the C$^*$-subalgebra of $A$
generated by $p$ and $q$ is isometrically to  $\CC \oplus^{\infty}
\CC$ when $p$ and $q$ are orthogonal and isomorphic to $M_2
(\CC)$, in other case.  More concretely, by \cite[Theorem 1.3
]{RaSin} (see also \cite[\S 3]{Ped68}),  denoting by $C_{p,q}$
the C$^*$-subalgebra of $A$ generated by $p$ and $q$, one of the
following statements holds: \begin{enumerate}[$a)$] \item $p\perp
q$ and there exists an isometric C$^*$-isomorphism $\Phi : C_{p,q}
\to \CC \oplus^{\infty} \CC$ such that $\Phi (p) = (1,0)$ and
$\Phi (q) = (0,1)$. \item  $p$ and $ q$ are not orthogonal and
there exist $0<t<1$ and an isometric C$^*$-isomorphism $\Phi :
C_{p,q} \to M_2 (\CC)$ such that $\Phi (p) = \left(
\begin{array}{cc}
                                                     1 & 0 \\
                                                     0 & 0 \\
                                                   \end{array}
                                                 \right)$ and  $\Phi (q) = \left(
                                                                             \begin{array}{cc}
                                                                               t & \sqrt{t(1-t)} \\
                                                                               \sqrt{t(1-t)} & 1-t \\
                                                                             \end{array}
                                                                           \right).$
\end{enumerate}\smallskip

In the setting of JB$^*$-algebras we have:

\begin{lem} \label{jbsubalgebra}Let $p$ and $q$ be two minimal projections in a
JB$^*$-algebra $J$ with $q\neq p$ and let $J_{p,q}$ denote the
JB$^*$-subalgebra of $J$ generated by $p$ and $q.$ Then one of the following statements holds:
\begin{enumerate}[$a)$] \item $p\perp q$ and there exists an isometric JB$^*$-isomorphism
$\Phi : J_{p,q} \to \CC \oplus^{\infty} \CC$ such that $\Phi (p) = (1,0)$ and $\Phi (q) = (0,1)$.
\item  $p$ and $ q$ are not orthogonal and there exist $0<t<1$ and an isometric JB$^*$-isomorphism
$\Phi : C \to S_2 (\CC)$ such that $\Phi (p) = \left( \begin{array}{cc}
                                                     1 & 0 \\
                                                     0 & 0 \\
                                                   \end{array}
                                                 \right)$
and  $\Phi (q) = \left( \begin{array}{cc}
  t & \sqrt{t(1-t)} \\
   \sqrt{t(1-t)} & 1-t \\
   \end{array}
   \right),$ where $S_2 (\CC)$ denotes the type 3 Cartan factor of all
   symmetric operators on a two dimensional complex Hilbert space.
\end{enumerate}

\noindent Moreover, the JB$^*$-subtriple of $J$ generated by $p$
and $q$ coincides with $J_{p,q}$.
\end{lem}

\begin{proof} Statement $a)$ is clear.  Let us assume that $p$ and $q$ are not orthogonal. The Shirshov-Cohn theorem (see \cite[Theorem 7.2.5]{HoS}) assures that $J_{p,q}$ is a JC$^*$-algebra, that is,  a Jordan *-subalgebra of some C$^*$-algebra $A$. The symbol $C_{p,q}$ will stand for the (associative) C$^*$-subalgebra of $A$ generated by $p$ and $q$. We have already commented that  there exist $0<t<1$ and an isometric C$^*$-isomorphism $\Phi : C_{p,q} \to M_2 (\CC)$ such that $\Phi (p) = \left( \begin{array}{cc}
                                                     1 & 0 \\
                                                     0 & 0 \\
                                                   \end{array}
                                                 \right)$ and  $\Phi (q) = \left(
                                                                             \begin{array}{cc}
                                                                               t & \sqrt{t(1-t)} \\
                                                                               \sqrt{t(1-t)} & 1-t \\
                                                                             \end{array}
                                                                           \right).$\smallskip

Since $J_{p,q}$ is a Jordan *-subalgebra of $C_{p,q}$, $J_{p,q}$ can be identified with the Jordan *-subalgebra of $M_2 (\CC)$ generated by the matrices $$P:= \left( \begin{array}{cc}
                                                     1 & 0 \\
                                                     0 & 0 \\
                                                   \end{array}
                                                 \right) \hbox{ and } Q:=\left(
                                                                             \begin{array}{cc}
                                                                               t & \sqrt{t(1-t)} \\
                                                                               \sqrt{t(1-t)} & 1-t \\
                                                                             \end{array}
                                                                           \right).$$\smallskip

It can be easily checked that $P\circ Q =\left(
                                                                             \begin{array}{cc}
                                                                               t & \frac12 \sqrt{t(1-t)} \\
                                                                              \frac12 \sqrt{t(1-t)} & 0 \\
                                                                             \end{array}
                                                                           \right),$ $2 P\circ Q - 2t P= \left(
                                                                             \begin{array}{cc}
                                                                               0 & \sqrt{t(1-t)} \\
                                                                               \sqrt{t(1-t)} & 0\\
                                                                             \end{array}
                                                                           \right),$   and $Q- (2 P\circ Q - 2t P) - t P  =\left(
                                                                             \begin{array}{cc}
                                                                               0 & 0 \\
                                                                               0 & 1-t \\
                                                                             \end{array}
                                                                           \right).$
These identities show that  $J_{p,q}$ contains the generators of
the JB$^*$-algebra $S_2 (\CC)$, and hence $J_{p,q}$ identifies
with $S_2 (\CC)$.\smallskip

Ir order to prove the last assertion, let $E_{p,q}$ denote the JB$^*$-subtriple of $J$
generated by $p$ and $q.$ As $J_{p,q}$ is itself a subtriple containing $p$ and $q$, we have
$E_{p,q}\subseteq J_{p,q}.$ If $p\perp q$ then it can easily be seen that
$E_{p,q}\cong \CC \oplus^\infty \CC \cong J_{p,q} .$ Let us assume that $p$
and $q$ are not orthogonal.\smallskip

From Proposition 5 in \cite{FriRu}, $E_{p,q}$ is isometrically
JB$^*$-triple isomorphic to $M_{1,2} (\CC)$ or $S_2(\CC)$. If
$E_{p,q}$ is a rank-one JB$^*$-triple, that is, $E\cong M_{1,2}
(\CC)$, then $P_0 (p) (q) $ must be zero. Thus, according to the
above representation, we have $1-t=0$, which is impossible.
\end{proof}

A JB$^*$-algebra which is a weakly compact JB$^*$-triple will be
called \emph{weakly compact} or \emph{dual} (see \cite{Bu82}).
Every positive element $x$ in a weakly compact JB$^*$-algebra $J$
can be written in the form $x= \sum_n \lambda_n p_n$
for a suitable $(\lambda_n) \in c_0$ and a family of mutually orthogonal minimal
projections $(p_n)$ in $J$ (see Theorem 3.3 in \cite{Bu82}).\smallskip

Our next theorem extends \cite[Theorem 11]{BurGarPe}.

\begin{thm}
\label{t authom cont dual j-alg} Let $T: J\rightarrow E$ be a
biorthogonality preserving linear surjection from a weakly compact
JB$^*$-algebra onto a JB$^*$-triple. Then $T$ is continuous and
$\|T\|\leq 2\sup \left\{ \|T(p)\| : p \hbox{ minimal projection in
} J\right\}$.
\end{thm}

\begin{proof} Since $J$ is a JB$^*$-algebra, it is enough to show
that $T$ is bounded on positive norm-one elements. In this case,
it suffices to prove that the set
$$\mathcal{P} = \left\{ \|T(p)\| : p \hbox{ minimal projection in
} J\right\} $$ is bounded (compare the proof of Proposition \ref{t
weakly compact}).\smallskip

Suppose, on the contrary, that $\mathcal{P}$ is unbounded. We
shall show by induction that there exists a sequence $(p_n)$ of
mutually orthogonal minimal projections in $J$ such that
$\|T(p_n)\| > n$.\smallskip

The case $n=1$ is clear. The induction hypothesis guarantees the
existence of mutually orthogonal minimal projections
$p_1,\ldots,p_n$ in $J$ with $\|T(p_k)\| >k$ for all $k\in
\{1,\ldots,n\}$.\smallskip

By assumption, there exists a minimal projection $q\in J$
satisfying $$\|T(q) \| > \max \{
\|T(p_1)\|,\ldots,\|T(p_n)\|,n+1\}.$$ We claim that $q$ must be
orthogonal to each $p_j$. If that is not the case, there exists
$j$ such that $p_j$ and $q$ are not orthogonal. Let $C$ denote the
JB$^*$-subtriple of $J$ generated by $q$ and $p_j$. We conclude from
Lemma \ref{jbsubalgebra} that $C$ is isomorphic to the
JB$^*$-algebra $S_2(\CC).$\smallskip

Let  $g_1=\left(
\begin{array}{cc}
1 & 0 \\
0 & 0 \\
\end{array}
\right)$ and $g_2= \left(
\begin{array}{cc}
0 & 0 \\
0 & 1 \\
\end{array}
\right)$. Then $g_1+g_2$ is the unit element in $C\cong S_2(\CC)$.
By Theorem \ref{thm minimal.to.minimal},
$w_1:=\frac{1}{\|T(g_1)\|} T(g_1)$ and $w_2:=\frac{1}{\|T(g_2)\|}
T(g_2)$ are two orthogonal minimal tripotents in $E$. The element
$w= w_1 +w_2$ is a rank-2 tripotent in $E$ and coincides with the
range tripotent of the element $h=T(g_1+g_2)= \|T(g_1)\| w_1 + \|T(g_2)\| w_2$.
Furthermore, $h$ is invertible in $E_2 (w)$, and
by Theorem \ref{t BurFerGarPe} (see also \cite[Corollary 4.1 $(b)$]{BurFerGarPe}),
${T(C)} \subseteq E_{2} (w)$.

The rest of the argument is  parallel to the argument in the proof of
Theorem \ref{biopcartan}.\smallskip

The finite dimensionality of the JB$^*$-subtriple $C$ assures that $T(C)$ is norm closed and
$T|_{C}: C\cong S_2(\CC) \to E$ is a continuous biorthogonality preserving linear operator.
Theorem \ref{t BurFerGarPe} guarantees the existence of a Jordan $^*$-homomorphism $S: C \to
E_2 (w) $ satisfying that $S (g_1+g_2) = w$, $S(C)$ and $h$
operator commute and \begin{equation}\label{eq new b} T(z) = h\circ_w S(z), \hbox{ for all } z\in C.
\end{equation} It follows from the operator commutativity of $h^{-1}$ and $S(C)$ that $S(z) = h^{-1}\circ_w T(z)$
for all $z\in C$. The injectivity of $T$ implies that $S$ is a Jordan *-monomorphism.\smallskip

Lemma 2.7 in \cite{FerMarPe2}
shows that $E_2 (w) = E_2 (w_1 +w_2)$ coincides with
$\mathbb{C}\oplus^{\ell_{\infty}}\mathbb{C}$ or with a spin
factor. Since $3= \dim ({T(C)}) \leq \dim(E_{2} (w))$,
we deduce that $E_2 (w)$ is a spin factor with inner product
$(.|.)$ and conjugation $x\mapsto \overline{x}$. We may assume, by Remark
\ref{r min trips in spin}, that $(w_1|w_1) = \frac12,$ $(w_1|\overline{w}_1)=0$,
and $w_2 = \overline{w}_1$.\smallskip

Now, taking $g_3=\left(
\begin{array}{cc}
0 & 1 \\
1 & 0 \\
\end{array}
\right) \in C\cong S_2(\CC)$, the element $w_3:= S(g_3)$ is a tripotent
in $E_2 (w)$ with $\J{w_i}{w_i}{w_3} = \frac12 w_3,$ for every $i \in
\{1,2\}$. Remark \ref{r min trips in
spin} implies that $(w_3|w_1) =(w_3|w_2)=0$. Let $M$ denote the JB$^*$-subtriple of
$E_2 (w)$ generated by $w_1,w_2,$ and $w_3$. The mapping $S: C\cong S_2(\CC) \to M$ is
a Jordan *-isomorphism.\smallskip

Combining $(\ref{eq new b})$ and $(\ref{eq spin product})$ we get
$$T(g_3) = h\circ_{w} S(g_3) = \J hw{w_3} = \frac{\|T(g_1)\|+\|T(g_2)\|}{2} w_3.$$\smallskip
Since $T(g_1) = \|T(g_1)\|\ w_1$, $T(g_2) = \|T(g_2)\|\ w_2$, and $C$ is linearly
generated by $g_1,g_2$ and $g_3,$ we deduce that $T(C) \subseteq M$ with
$3= \dim(T(C)) \leq  \dim(M) = 3$. Thus, $T(C) = M$ is a JB$^*$-subtriple of $E$.\smallskip

The mapping $T|_{C} : C\cong S_2(\CC) \to T(C)$ is a continuous
biorthogonality preserving linear bijection.
Theorem \ref{charbbiops} guarantees the existence of
a scalar $\lambda\in \CC\backslash\{0\}$ and a triple isomorphism
$\Psi : C \to T(C)$ such that $T(x) = \lambda \Psi (x)$ for all
$x\in C$. Since $p_j$ and $q$ are projections, $\|\Psi (q)\| = \|\Psi (p_j)\|=1$.
Hence $\|T(p_j) \| = |\lambda|$ and $\|T(q) \| = |\lambda|,$ contradicting the induction hypothesis.
Therefore $q\perp p_j$, for every $j= 1,\ldots,n$.\smallskip

It follows by induction that there exists a sequence  $(p_n)$ of
mutually orthogonal minimal projections in $J$ such that
$\|T(p_n)\| > n$. The series $\sum_{n=1}^{\infty}
\frac{1}{\sqrt{n}} p_n$ defines an element $a$ in $J$ (compare
Remark \ref{r Cauchy series}). For each natural $m$, $a$
decomposes as the orthogonal sum of $\frac{1}{\sqrt{m}} p_m$ and
$\sum_{n\neq m}^{\infty} \frac{1}{\sqrt{n}} p_n$, therefore
$$T(a)= \frac{1}{\sqrt{m}} T(p_m) + T\left(\sum_{n\neq m}^{\infty}
\frac{1}{\sqrt{n}} p_n\right),$$ with $\frac{1}{\sqrt{m}} T(p_m)
\perp T\left(\sum_{n\neq m}^{\infty} \frac{1}{\sqrt{n}}
p_n\right).$ This argument implies that $$\|T(a)\|=\max \left\{
\frac{1}{\sqrt{m}} \left\|T(p_m) \right\| ,\
\left\|T\left(\sum_{n\neq m}^{\infty} \frac{1}{\sqrt{n}}
p_n\right)\right\| \right\}> \sqrt{m}.$$ Since $m$ was arbitrarily
chosen, we have arrived at the desired contradiction.
\end{proof}

By Proposition 2 in \cite{HoMarPeRu}, every Cartan factor of type
1 with dim$(H) = \hbox{dim} (H')$, every Cartan factor of type 2
with dim$(H)$ even, or infinite, and every Cartan factor of type 3
is a JBW$^*$-algebra factor for a suitable Jordan product and
involution. In the case of $C$ being a Cartan factor which is also
a JBW$^*$-algebra, then the corresponding elementary JB$^*$-triple
$K(C)$ is a weakly compact JB$^*$-algebra.

\begin{cor}
\label{c authom cont dual j-alg factor} Let $K$ be an elementary
JB$^*$-triple of type 1 with dim$(H) = \hbox{dim} (H')$, or of
type 2 with dim$(H)$ even, or infinite, or of type 3. Suppose that
$T: K\rightarrow E$ is a biorthogonality preserving linear
surjection from $K$ onto a JB$^*$-triple. Then $T$ is continuous.
Further, since $K^{**}$ is a JBW$^{*}$-algebra factor, Theorem
\ref{charbbiops} assures that $T$ is a scalar multiple of a triple
isomorphism.$\hfill\Box$
\end{cor}

\begin{thm}\label{tautomwcomp} Let $T: E\rightarrow F$ be a
biorthogonality preserving linear surjection between two
JB$^*$-triples, where $E$ is weakly compact containing no infinite
dimensional rank-one summands. Then $T$ is continuous.
\end{thm}

\begin{proof}
Since $E$ is a weakly compact JB$^*$-triple, the statement follows
from Proposition \ref{t weakly compact} as soon as we prove that
the set $$\mathcal{T}:=\left\{ \|T(e)\| : e \hbox{ minimal
tripotent in } E\right\} $$ is bounded.\smallskip

We know that $E=\bigoplus^{c_0}_{\alpha \in \Gamma} K_{\alpha}$,
where $\{K_\alpha: \alpha\in \Gamma \} $ is a family of elementary
JB$^*$-triples (see Lemma 3.3 in \cite{BuChu}). Now, Lemma \ref{l
basic prop annihilator} guarantees  that $T(K_\alpha) = T(K_\alpha^{\perp\perp}) =
T(K_\alpha)^{\perp\perp}$ is a norm closed inner ideal, for every $\alpha\in \Gamma$.\smallskip

For each $\alpha\in \Gamma$, $K_\alpha$ is finite dimensional, or a type 1 elementary JB$^*$-triple of rank bigger than one, or a JB$^*$-algebra. It follows, by Corollary
\ref{scalartripletypeI} and Theorem \ref{t authom cont dual j-alg}, that $T|_{K_{\alpha}} : K_{\alpha} \to T(K_{\alpha})$ is continuous. 
\smallskip

Suppose that $\mathcal{T}$ is unbounded.  Having in mind that every minimal tripotent in $E$ belongs to a unique factor $K_\alpha$, by Proposition \ref{t weakly compact}, there exists a sequence $(e_n)$ of mutually orthogonal minimal tripotents in $E$ such that $\|T(e_n)\|$ diverges to $+\infty$. The element $z:= \sum_{n=1}^{\infty} {\|T(e_n)\|}^{-\frac12} \  {e_n}$ lies in $E$ and hence $\|T(z)\| <\infty$. We fix an arbitrary natural $m$. Since $z-{\|T(e_m)\|}^{-\frac12} \  {e_m}$  and ${\|T(e_m)\|}^{-\frac12} \  {e_m}$ are orthogonal, we have $$T(z-{\|T(e_m)\|}^{-\frac12} \  {e_m}) \perp T( {\|T(e_m)\|}^{-\frac12} \  {e_m}),$$ and hence $$\|T(z) \| = \|T(z -{\|T(e_m)\|}^{-\frac12} \  {e_m})) + T({\|T(e_m)\|}^{-\frac12} \  {e_m})\| $$ $$= \max \left\{\|T(z-{\|T(e_m)\|}^{-\frac12} \  {e_m})\|, {\|T(e_m)\|}^{-\frac12} \   \|T({e_m})\|\right\} \geq \sqrt{\|T(e_m)\|},$$ which contradicts that ${\|T(e_m)\|}^{\frac12} \to +\infty$. Therefore $\mathcal{T}$ is bounded.
\end{proof}

\begin{cor}
\label{c authom cont K(E)} Let $T: E\rightarrow F$ be a
biorthogonality preserving linear surjection between two
JB$^*$-triples, where $K(E)$ contains no infinite
dimensional rank-one summands. Then $T|_{K(E)}: K(E) \to K(F)$ is continuous.
\end{cor}

\begin{proof} Let us take $x\in K(E)$. In this case $x$ can be written in
the form $x = \sum_{n} \lambda_n u_n,$ where $u_n$ are mutually
orthogonal minimal tripotents of $E$, and $\|x\| = \sup
\{|\lambda_n| : n\}$ (compare Remark 4.6 in \cite{BuChu}). For
each natural $m$ we define $y_m := T \left( \sum_{n\geq
m+1}^{\infty} \lambda_n u_n \right)$. Theorem \ref{thm
minimal.to.minimal} guarantees that $T(x_m) =
T\left(\sum_{n=1}^{m} \lambda_n u_n\right)$ defines a sequence in
$K(F)$.\smallskip

Since, by Lemma \ref{l lim series to zero}, $y_m \to 0$ in norm,
we deduce that $T(x_m) = T(x) - y_m$ tends to $T(x)$ in norm.
Therefore $T(K(E)) = K(F)$ and $T|_{K(E)}: K(E) \to K(F)$ is a
biorthogonality preserving linear surjection between two weakly
compact JB$^*$-triples. The result follows now from Theorem
\ref{tautomwcomp} above.
\end{proof}

\begin{rem}\label{countercontinweakly}{\rm Remark 15 in \cite{BurFerGarMarPe} already
pointed that  the conclusion of Theorem \ref{tautomwcomp} is no
longer true if we  allow $E$ to have infinite dimensional rank-one
summands. Indeed, let $E=L(H)\oplus^\infty L(H,\CC)$, where $H$ is
an infinite dimensional complex Hilbert space. We can always find
an unbounded bijection $S: L(H,\CC)\rightarrow L(H,\CC).$ Since
$L(H,\CC)$ is a rank-one JB$^*$-triple then $S$ is a biorthogonality
preserving linear bijection and the mapping $T:E \rightarrow E,$
given by $x+y\mapsto x+T(y),$ satisfies the same properties.}
\end{rem}

\begin{cor} \label{isomweaklycom} Two weakly compact JB$^*$-triples
containing no rank-one summands are isomorphic if, and only if,
there exists a biorthogonality preserving linear surjection between
them.
\end{cor}

\section{Biorthogonality preservers between atomic JBW$^*$-triples}

A JBW$^*$-triple $E$ is said to be \emph{atomic} if it coincides
with the weak$^*$ closed ideal generated by its minimal
tripotents. Every atomic  JBW$^*$-triple can be written as a
$l_\infty$-sum of Cartan factors \cite{FriRu86}. \smallskip

The aim of this section is to study when the existence of a
biorthogonality preserving linear surjection between two atomic
JBW$^*$-triples implies that they are isomorphic (note that
continuity is not assumed). We shall establish an automatic
continuity result for biorthogonality linear surjections between
atomic JBW$^*$-triples containing no rank-one factors.\smallskip

Before dealing with the main result, we survey some results
describing the elements in the predual of a Cartan factor. We make
use of the description of the predual of $L(H)$ in terms of the
\emph{trace class} operators (compare \cite[\S II.1]{Tak}). The
results, included here for completeness reasons, are direct
consequences of this description but we do not know a explicit
reference.\smallskip

Let $C=L(H,H')$ be a type 1 Cartan factor. Lemma 2.6 in \cite{Pe}
assures that  each $\varphi$ in $C_*$ can be written in the form
$\varphi := \sum_{n=1}^{\infty} \lambda_n \varphi_n$, where
$(\lambda_n)$ is a sequence in $\ell_1^{+}$ and each  $\varphi_n$
is an extreme point of the closed unit ball of $C_*$. More
concretely, for each natural $n$ there exist norm-one elements
$h_n\in H$ and $k_n\in H'$ such that $\varphi_n (x) = (x(h_n) |
k_n)$, for every $x\in C$, that is, for each natural $n$ there
exists a minimal tripotent $e_n$ in $C$ such that $P_2 (e_n) (x) =
\varphi_n (x) e_n$, for every $x\in C$ (compare \cite[Proposition
4]{FriRu}).\smallskip

We consider now (infinite dimensional) type 2 and type 3 Cartan
factors. Let  $j$ be a conjugation on a complex Hilbert space $H$,
and let ``$t$'' denote the linear involution on $L(H)$ defined by
$x\mapsto x^t:=jx^*j$. Let $C_2 = \{x\in L(H) : x^t =-x\}$ and
$C_3 = \{x\in L(H) : x^t =x\}$ be Cartan factors of type 2 and 3,
respectively.\smallskip

Noticing that $L(H) = C_2 \oplus C_3$, it is easy to see that
every element $\varphi$ in $(C_2)_*$ (respectively, $(C_3)_*$)
admits an extension of the form $\widetilde{\varphi} = \varphi
\pi$, where $\pi$ denotes the canonical projection of $L(H)$ onto
$C_2$ (respectively, $C_3$). Making use of \cite[Lemma 1.5]{Tak},
we can find an element $x_{_{\widetilde{\varphi}}}\in K(H)$
satisfying
\begin{equation}
\label{eq trace class}(x_{_{\widetilde{\varphi}}} (h) | k)=
\widetilde{\varphi} (h\otimes k), \ (h,k\in H).
\end{equation}

Since, for each $x\in L(H)$, $\widetilde{\varphi} (x) = \frac12 \
\widetilde{\varphi} (x-x^t) $, we can easily check, via $(\ref{eq
trace class})$, that $x_{_{\widetilde{\varphi}}}^t = -
x_{_{\widetilde{\varphi}}}$. Therefore
$x_{_{\widetilde{\varphi}}}\in K_2 = K(C_2)$. From \cite[Remark
4.6]{BuChu} it may be concluded that $x_{_{\widetilde{\varphi}}} $
can be (uniquely) written as a norm convergent (possibly finite)
sum $x_{_{\widetilde{\varphi}}} = \sum_{n} \lambda_n u_n$, where
$u_n$ are mutually orthogonal minimal tripotents in $K_2$ and
$(\lambda_n) \in c_0$ (It should be noticed here that $u_n$ is a
minimal tripotent in $C_2$ but it needs not be minimal in $L(H)$.
In any case, $u_n$ is minimal in $L(H)$ or it can be written as a
convex combination of two minimal tripotents in $L(H)$). For each
$(\beta_n)\in c_0$, $z:= \sum_{n} \beta_n u_n \in K_2$ and, by
$(\ref{eq trace class})$, $\sum_{n} \lambda_n \beta_n =
\widetilde{\varphi} (z)= \varphi (z) < \infty$. Thus, $(\lambda_n)
\in \ell_1,$ and, another application of $(\ref{eq trace class})$
shows that $\varphi (x) = \sum_{n} \lambda_n \varphi_n (x)$,
$(\forall x\in C_2),$ where $\varphi_n$ lies in $(C_2)_*$ and
satisfies $P_2 (u_n) (x) = \varphi_n (x) u_n$. A similar reasoning
remains true for $C_3$.\smallskip

We have therefore prove:

\begin{proposition}\label{p trace class types 1,2,3} Let $C$ be an
infinite dimensional Cartan factor of type 1, 2 or 3. For each
$\varphi$ in $C_*$, there exist a sequence $(\lambda_n) \in
\ell_1,$ and a sequence $(u_n)$ of mutually orthogonal minimal
tripotents in $C$ such that $$\|\varphi\| = \sum_{n=1}^{\infty}
|\lambda_n| \hbox{ and }\varphi (x) = \sum_n \lambda_n \ \varphi_n
(x), \ (x\in C),$$ where for each $n\in \NN$, $\varphi_n (x)\ u_n
= P_2 (u_n) (x), \ (x\in C)$. $\hfill\Box$
\end{proposition}

Let $T: E\to F$ be a biorthogonality preserving linear surjection
between atomic JBW$^*$-triples, where $E$ contains no rank-one
Cartan factors. In this case $K(E)$ and $K(F)$ are weakly compact
JB$^{*}$-triples with $K(E)^{**} = E$ and  $K(F)^{**} = F$.
Corollary \ref{c authom cont K(E)} assures that $T|_{K(E)}: K(E)
\to K(F)$ is continuous. This is not, a priori, enough reason to
guarantee that $T$ is continuous. In fact, for each non reflexive
Banach space $X$ there exists an unbounded linear operator $S:
X^{**}\to X^{**}$ satisfying that $S|_{X} : X\to X$ is continuous.
The main result of this section establishes that a mapping $T$ in
the above conditions is automatically continuous.

\begin{thm}
\label{t atomic} Let $T: E\to F$ be a biorthogonality preserving
linear surjection between atomic JBW$^*$-triples,
where $E$ contains no rank-one Cartan factors. Then $T$ is continuous.
\end{thm}

\begin{proof} Corollary \ref{c authom cont K(E)}
assures that $T|_{K(E)}: K(E) \to K(F)$ is continuous. By Lemma
3.3 in \cite{BuChu} $K(E)$ decomposes as a $c_{0}$-sum of all
elementary triple ideals of $E$, that is, if $E=
\oplus^{\ell_{\infty}} C_{\alpha},$ where each $C_{\alpha}$ is a
Cartan factor, then $K(E) = \oplus^{c_{0}} K(C_{\alpha}).$ By
Proposition \ref{p direct sums}, for each $\alpha$,
$T(K_{\alpha})$ (respectively, $T(C_{\alpha}$) is a norm closed
(respectively, weak$^{*}$ closed) inner ideal of $K(F)$
(respectively, $F$) and $K(F) = \oplus^{c_{0}} T(K(C_{\alpha}))$
(respectively,  $F = \oplus^{c_{0}} T(C_{\alpha})$).\smallskip

For each $\alpha$, $C_{\alpha}$ is finite dimensional or an
infinite dimensional Cartan factor of type 1, 2 or 3. Corollaries
\ref{scalartripletypeI} and \ref{c authom cont dual j-alg factor}
prove that  the operator $T|_{{K(C_{\alpha})}} : K(C_{\alpha}) \to
T(K(C_{\alpha}))$ is a scalar multiple of a triple isomorphism. We
claim that, for each $\alpha$ and each $\varphi_{\alpha}$ in the
predual of $T(C_{\alpha})$, $\varphi_{\alpha} T$ is weak$^*$
continuous. There is no lost of generality in assuming that
$C_{\alpha}$ is infinite dimensional.\smallskip

Each minimal tripotent $f$ in $F$, lies in a unique elementary
JB$^{*}$-triple $T(K(C_{\alpha}))$. Since $T|_{{K(C_{\alpha})}} :
K(C_{\alpha}) \to T(K(C_{\alpha}))$ is a scalar multiple of a
triple isomorphism, there exist a non-zero scalar
$\lambda_{\alpha}$ and a minimal tripotent $e$ satisfying $ T^{-1}
(f) = \lambda_{\alpha} e $, $|\lambda_{\alpha}|\leq
\|\left(T|_{{K(C_{\alpha})}}\right)^{-1}\|\leq
\|\left(T|_{{K(E)}}\right)^{-1}\|,$ and
\begin{equation} \label{eq 0,2 peirce weakly compact preserved}T
(K(C_{\alpha})_{i} (e) ) =T(K(C_{\alpha}))_i (f),
\end{equation} for every
$i=0,1,2$.  Theorem \ref{thm minimal.to.minimal}  shows that $T
((C_{\alpha})_{i} (e) ) =T(C_{\alpha})_i (f),$ for every $i=0,2$.
Since $K(E)$ is an ideal of $E$ and $e$ is a minimal tripotent,
$(C_{\alpha})_{1} (e) = E_{1} (e) =  K(E)_{1} (e) =
K(C_{\alpha})_{1} (e).$ It follows by $(\ref{eq 0,2 peirce weakly
compact preserved})$ that $$T ((C_{\alpha})_{i} (e) )
=T((C_{\alpha}))_i (f),$$ for every $i=0,1,2$. Consequently, $P_2
(f) T  = \lambda_{\alpha}^{-1} P_2 (e) \in (C_{\alpha})_{*}$, and
$|\lambda_{\alpha}^{-1}|\leq \|T|_{{K(C_{\alpha})}}\|\leq
\|T|_{{K(E)}}\|.$ \smallskip

Since $f$ was arbitrarily chosen among the minimal tripotents in
$F$, (equivalently, in $T(K(C_{\alpha}))$), Proposition \ref{p
trace class types 1,2,3} assures that $\varphi_{\alpha} T \in
E_*$, with $\|\varphi_{\alpha} T\| \leq \|T|_{{K(E)}}\|$, for
every $\varphi_{\alpha}\in \left(T(C_{\alpha})\right)_{*}.$
Therefore, $T$ is bounded with $$\|T\| \leq \|T|_{{K(E)}}\| \leq
\|T\|.$$
\end{proof}

\begin{cor} \label{c isom atomic} Two atomic JBW$^*$-triples
containing no rank-one summands are isomorphic if, and only if,
there is a biorthogonality preserving linear surjection between
them.$\hfill\Box$
\end{cor}

The conclusion of the above Theorem does not hold for atomic
JBW$^*$-triples containing rank-one summands.

\subsection*{Acknowledgements}
The authors would like to express their gratitude to the anonymous
referee whose valuable comments improved the final form of this paper.\smallskip

This research was partly supported by I+D MEC project 
no. MTM 2008-02186, and Junta de Andaluc{\'\i}a 
grants FQM 0199 and FQM 3737.

\end{document}